\documentclass[a4paper,12pt]{amsart}
\usepackage{amssymb}
\usepackage{ifthen}
\newtheorem{thm}{Theorem}
\newtheorem{cor}{Corollary}
\newtheorem{lem}{Lemma}

\newtheorem{conj}{Conjecture}

\theoremstyle{definition}

\newtheorem{prob}[equation]{Problem}

\newenvironment{rem}{%
\bigskip
\noindent \textsl{{\sl Remark. }}}{\bigskip}
\newenvironment{rems}{%
\bigskip
\noindent \textsl{{\sl Remarks. }}}{\bigskip}

\newcounter {own}
\def\theown {\thesection       .\arabic{own}}

\newenvironment{pf}[1][]{%
 \vskip 3mm
 \noindent
 \ifthenelse{\equal{#1}{}}%
  {{\slshape Proof. }}%
  {{\slshape #1.} }%
 }%
{\qed\bigskip}

\newcounter{alphabet}
\newcounter{tmp}

\newcommand{\ID}{{\mathbb D}}

\newcommand{\IC}{{\mathbb C}}

\def\be{\begin{equation}}
\def\ee{\end{equation}}

\newcommand{\bee}{\begin{enumerate}}
\newcommand{\eee}{\end{enumerate}}

\newcommand{\pay}{\!\!\!}
\newcommand{\blem}{\begin{lem}}
\newcommand{\elem}{\end{lem}}
\newcommand{\bthm}{\begin{thm}}
\newcommand{\ethm}{\end{thm}}
\newcommand{\bcor}{\begin{cor}}
\newcommand{\ecor}{\end{cor}}
\newcommand{\beg}{\begin{examp}}
\newcommand{\eeg}{\end{examp}}
\newcommand{\begs}{\begin{examples}}
\newcommand{\eegs}{\end{examples}}
\newcommand{\bdefe}{\begin{defin}}
\newcommand{\edefe}{\end{defin}}
\newcommand{\bprob}{\begin{prob}}
\newcommand{\eprob}{\end{prob}}
\newcommand{\bei}{\begin{itemize}}
\newcommand{\eei}{\end{itemize}}

\newcommand{\bcon}{\begin{conj}}
\newcommand{\econ}{\end{conj}}
\newcommand{\bcons}{\begin{conjs}}
\newcommand{\econs}{\end{conjs}}
\newcommand{\bprop}{\begin{propo}}
\newcommand{\eprop}{\end{propo}}
\newcommand{\br}{\begin{rem}}
\newcommand{\er}{\end{rem}}
\newcommand{\brs}{\begin{rems}}
\newcommand{\ers}{\end{rems}}
\newcommand{\bo}{\begin{obser}}
\newcommand{\eo}{\end{obser}}
\newcommand{\bos}{\begin{obsers}}
\newcommand{\eos}{\end{obsers}}
\newcommand{\bpf}{\begin{pf}}
\newcommand{\epf}{\end{pf}}
\newcommand{\ba}{\begin{array}}
\newcommand{\ea}{\end{array}}
\newcommand{\beq}{\begin{eqnarray}}
\newcommand{\beqq}{\begin{eqnarray*}}
\newcommand{\eeq}{\end{eqnarray}}
\newcommand{\eeqq}{\end{eqnarray*}}

\newcommand{\ds}{\displaystyle}

\begin{document}
\bibliographystyle{amsplain}
\title{Domains of variability of Laurent coefficients and the convex hull for
the family of concave univalent functions}
\author{B. Bhowmik}
\address{B. Bhowmik, Department of Mathematics,
Indian Institute of Technology Madras, Chennai-600 036, India.}
\email{ditya@iitm.ac.in}
\author{S. Ponnusamy}
\address{S. Ponnusamy, Department of Mathematics,
Indian Institute of Technology Madras, Chennai-600 036, India.}
\email{samy@iitm.ac.in}
\author{K.-J. Wirths}
\address{K.-J. Wirths, Institut f\"ur Analysis und Algebra, TU Braunschweig,
38106 Braunschweig, Germany}
\email{kjwirths@tu-bs.de}


\begin{abstract}
Let $\ID$ denote the open unit disc and let $p\in (0,1)$. We consider the
family $Co(p)$ of functions $f:\ID\to \overline{\IC}$ that satisfy the
following conditions:
\bee
\item[(i)] $f$ is meromorphic in $\ID$ and has a simple pole at the point $p$.
\item[(ii)] $f(0)=f'(0)-1=0$.
\item[(iii)] $f$ maps $\ID$ conformally onto a set whose complement with respect to
$\overline{\IC}$ is convex.
\eee
We determine the exact domains of variability of some coefficients $a_n(f)$
of the Laurent expansion
$$f(z)=\sum_{n=-1}^{\infty} a_n(f)(z-p)^n,\quad |z-p|<1-p,
$$
for $f\in Co(p)$ and certain values of $p$.
Knowledge on these Laurent coefficients is used to disprove a conjecture
of the third author on the closed convex hull of $Co(p)$ for certain values of $p$.
\end{abstract}
\thanks{}

\maketitle
\pagestyle{myheadings}
\markboth{B. Bhowmik, S. Ponnusamy, K.-J. Wirths}{Concave univalent function}

\bigskip

Let $\ID$ denote the open disc and let $p\in (0,1)$. We consider the family $Co(p)$ of
 functions $f:\ID\to \overline{\IC}$ that satisfy the following conditions:
 \bee
\item[(i)] $f$ is meromorphic in $\ID$ and has a simple pole at the point $p$.
\item[(ii)] $f(0)=f'(0)-1=0$.
\item[(iii)] $f$ maps $\ID$ conformally onto a set whose complement with respect to
$\overline{\IC}$ is convex.
\eee
In \cite{W} the third author of the present article proved the following theorem.

\bigskip
\noindent
{\bf Theorem A.} {\em
Let $p\in (0,1), f\in Co(p)$,  and let
$$ f(z)=\sum_{n=-1}^{\infty} a_n(f)(z-p)^n,\quad |z-p|<1-p,
$$
be the Laurent expansion of $f$ at the point $p$. Then the domain of variability of the
residuum $a_{-1}(f)$  is determined by the inequality
\be\label{f1}
\left|a_{-1}(f)\,+\,\frac{p^2}{1-p^4}\right|\,\leq\,\frac{p^4}{1-p^4}.
\ee
Equality is attained in {\rm (\ref{f1})}  if and only if
\be\label{f2}
f(z)=\frac{z\,-\,\frac{p}{1+p^2}\left(1+e^{i\theta}\right)z^2}
{\left(1-\frac{z}{p}\right)(1-zp)},\quad z\in \ID,
\ee
for some $\theta\in [0,2\pi]$.
}

\bigskip
Theorem A follows without difficulty from the following representation
theorem proved by Avkhadiev and Wirths in \cite{AW}.

\bigskip
\noindent
{\bf Theorem B.} {\em
Let $p\in (0,1)$. For any $f\in Co(p)$ there exists a function
$\omega:\ID\to \overline{\ID}$  holomorphic in $\ID$  such that
\be\label{f3}
f(z)=\frac{z\,-\,\frac{p}{1+p^2}(1+\omega(z))z^2}{\left(1-\frac{z}{p}\right)(1-zp)},
\quad z\in \ID.
\ee
}

\bigskip

On one hand, the present article originated in discussions among the authors whether
it is possible to derive the domains of variability of Laurent coefficients
$a_n(f), n\geq 0,$ for $f\in Co(p)$ from Theorem B.
On the other hand, the third author hoped that it would be possible to prove that the family of functions
defined by (\ref{f3}) represents the closed convex hull of $Co(p)$ in the topology
of uniform convergence on compact subsets of $\ID$.

In the sequel, we will determine the above domains of variability for $n=0$ and $n=1$
for certain values of $p$. For the remaining values of $p$ we will use these considerations
and some results of Livingston in \cite{L} to show that the above mentioned hope was in vain.

Our first result is an application of Theorem A and Theorem 4 in \cite{L}.

\bthm\label{th1}
 Let $p\in (0,1)$  and $f\in Co(p)$.  Then
\be\label{f4}
\mbox{\rm{Re}}\,\,a_0(f)\,\geq\,-\,\frac{p}{(1-p^2)^2}.
\ee
Equality is attained in {\rm (\ref{f4})} if and only if
\be\label{f5}
f(z)=\frac{z}{\left(1-\frac{z}{p}\right)(1-zp)},\quad z\in \ID.
\ee
\ethm\bpf
In  \cite[Theorem 4]{L}, Livingston proved that for $f\in Co(p)$ the inequality
$$\left|p\,+\,\frac{a_0(f)(1-p^2)}{a_{-1}(f)}\right|\,\leq\,\frac{1+p^2}{p}
$$
is valid. Hence, for any $f\in Co(p)$ there exists a number $\tau\in\overline{\ID}$
such that
\be\label{f6}
a_0(f)=\frac{a_{-1}(f)}{1-p^2}\left(-p\,+\,\tau\,\frac{1+p^2}{p}\right).
\ee
To prove (\ref{f4}) we have to determine the minimal real part of the product at the
right side of (\ref{f6}), where $a_{-1}(f)$ varies in the disc described by (\ref{f1}).
To that end it is sufficient to consider the points
$\tau = e^{i\varphi}, \varphi\in [0,2\pi],$ and to compute the minimum of the quantity

\vspace{8pt}

$\displaystyle \frac{-\,p}{(1-p^4)(1-p^2)}\left((1+p^2)\cos\varphi\,
-\,p^2\right)\,$
\[-\frac{p^3}{(1-p^4)(1-p^2)}\left((1+p^2)^2\sin^2\varphi+
\left((1+p^2)\cos\varphi\,-\,p^2\right)^2\right)^{1/2},
\]
where $\varphi\in [0,2\pi]$. Letting $x=\cos\varphi\in [-1,1]$ in this
expression and differentiating with respect to $x$ reveals there is no local extremum in
the interval $(-1,1)$. Therefore, it is easy to see that the minimum is attained for
$\tau =1$ and $a_{-1}(f) = -p^2/(1-p^2).$ According to Theorem A, this residuum occurs
only for the function (\ref{f5}) and for this function equality is attained in (\ref{f4}).
This concludes the proof of Theorem \ref{th1}.
\epf

For poles near the origin much more can be proved.

\bthm\label{th2} Let $p\in (0,\sqrt{3}-1]$  and $f\in Co(p)$. Then
the domain of variability of $a_0(f)$ is determined by the
inequality
\be\label{f7}
\left|\frac{1-p^2}{p}a_0(f)\,+\,\frac{1-p^2+p^4}{1-p^4}\right|\,\leq\,\frac{p^2(2-p^2)}{1-p^4}.
\ee
Equality is attained in {\rm (\ref{f7})}  if and only if $f$
is one of the functions given in {\rm (\ref{f2})}.
\ethm\bpf
We
multiply (\ref{f3}) by the denominator of the right side and
expand both side in power series with expansion point at $p$.
In the resulting equation, letting
$$\omega(z)\,=\,\sum_{n=0}^{\infty}c_n(z-p)^n,\quad z\in \ID,
$$
and comparing the constant terms and the coefficients of $(z-p)$,
we get
\be\label{f8}
a_{-1}(f)\,=\,\frac{-p^2}{1-p^4}\,+\,\frac{p^4}{1-p^4}\,c_0
\ee
and
\be\label{f9} a_{-1}(f)\,-\,\frac{1-p^2}{p}\,a_0(f)\,
=\,\frac{1-p^2}{1+p^2}\,-\,\frac{p^2}{1+p^2}\,(2c_0+pc_1).
\ee
It may be mentioned at this place that (\ref{f8}) and the inequality
$|c_0|\leq 1$ immediately prove Theorem A.

Further, we derive from (\ref{f8}) and (\ref{f9})
together the representation
\be\label{f10}
\frac{1-p^2}{p}\,a_0(f)\,+\frac{1-p^2+p^4}{1-p^4}\,=\,\frac{2p^2-p^4}{1-p^4}\,c_0\,
+\frac{p^3}{1+p^2}\,c_1.
\ee
Using the inequalities
$$|c_0|\leq 1\quad\mbox{\rm{and}}\quad|c_1|\leq\frac{1-|c_0|^2}{1-p^2},
$$
we get from (\ref{f10}) the inequality
$$ \left|\frac{1-p^2}{p}\,a_0(f)\,+\frac{1-p^2+p^4}{1-p^4}\right|
\leq \frac{p^2}{1-p^4}\left((2-p^2)|c_0|+p(1-|c_0|^2)\right).
$$
The function
$$g(x)=(2-p^2)x+p(1-x^2)
$$
has its local maximum at $x_M(p)=(2-p^2)/2p$. Since $x_M(p)\geq 1$ for
$p\in (0,\sqrt{3}-1]$, we get that
$$ \max\{g(x)\,\mid\,x\in[0,1]\}\,=\,g(1)=2-p^2
$$
for those $p$. This proves the inequality (\ref{f7}) for $f\in Co(p)$.
Obviously, $|c_0|=1$ implies that the only functions $f\in Co(p)$,
for which equality can occur there, are the functions (\ref{f2}).

The points in the disc described by (\ref{f7}) are attained for the
functions (\ref{f3}) with $\omega(z)\equiv c_0, |c_0|\leq 1$.
The fact that they belong to the class  $Co(p)$ has been proved in \cite{APW} and \cite{W}.
The proof of Theorem \ref{th2} is finished.
\epf

Now, we turn to the values of $p$ in the interval $(\sqrt{3}-1,1)$ and for them we get

\bthm\label{th3}
Let $p\in(\sqrt{3}-1,1)$. Then the closed convex hull of the class $Co(p)$
is a proper subset of the class of functions defined by {\rm (\ref{f3})}.
\ethm\bpf
It is a direct consequence of Theorem \ref{th1} that the coefficients
$a_0(f)$ of the functions
in the closed convex hull of $Co(p)$ satisfy the inequality (\ref{f4}), likewise.

On the other hand, let us insert into (\ref{f3}) the functions
\be\label{f11}
\omega_x(z)\,=\,\frac{-\left(\frac{z-p}{1-pz}\right)\,-\,x}
{1\,+\,x\left(\frac{z-p}{1-pz}\right)},\quad z\in \ID,
\ee
$x\in(0,1)$ fixed. A computation of the coefficients $a_0(f)$ for the resulting functions
using (\ref{f10}) delivers
$$ a_0(f)\,=\,\frac{-p}{(1-p^2)^2}\left(1\,+\,\frac{(1-x)p^2}{1+p^2}
\left(p(1+x)-(2-p^2)\right)\right).
$$
The right side is less than $-p/(1-p^2)^2$ for $x>(2-p^2-p)/p$ and
$(2-p^2-p)/p<1$ for $p\in(\sqrt{3}-1,1)$. Hence, the functions $f$ got by inserting
(\ref{f11}) into (\ref{f3}) do not belong to the closed convex hull of $Co(p)$ for
the values of $p$ indicated in Theorem \ref{th3} and $x\in ((2-p^2-p)/p,1)$.
\epf

In the sequel, we shall prove similar theorems as above concerning the coefficient $a_1(f)$.
During this program Theorem \ref{th1} may be replaced by the following theorem.

\bigskip
\noindent
{\bf Theorem C.} (see \cite[Theorem 3]{L}) {\em
Let $p\in (0,1)$  and $f\in Co(p).$
Then the inequality
$$ \left|a_1(f)\right|\,\leq\,\frac{p^2}{(1-p^2)^3}
$$
is valid.
}

\bigskip
Concerning the analogue to Theorem \ref{th2}, much more effort than before is needed because
of the appearance of $c_0, c_1,$ and $c_2$ in the formulas. To get a sharp result
nevertheless, we apply the  theory of extremum problems for
linear functionals on $H^p, 1\leq p \leq \infty,$ due to  Macintyre, Rogosinski,
and  Shapiro (see \cite{MR}, \cite{RS}, and Duren's book \cite{D} on $H^p$ spaces, Ch. 8).
This discussion enables us to prove

\bthm\label{th4}
 Let $p\in \left(0, 1-\frac{\sqrt{2}}{2}\right]$ and $f\in Co(p)$.
Then the domain of variability of $a_1(f)$ is determined by the inequality
\be\label{f12}
\left|a_1(f)\left(\frac{1-p^2}{p}\right)^2\,+\,\frac{p^2}{1-p^4}\right|\,
\leq\,\frac{1}{1-p^4}.
\ee
Equality is attained in {\rm (\ref{f12})}  if and only if $f$ is one of the
functions given in {\rm (\ref{f2})}.
\ethm\bpf
By the same procedure as in the proof of Theorem \ref{th2} we get in addition to
(\ref{f8}) and (\ref{f9}) comparing the coefficients of $(z-p)^2$
$$
a_0(f)-\,\frac{1-p^2}{p}a_1(f)\,=\,-\,\frac{p}{1+p^2}(1+c_0+2pc_1+p^2c_2).
$$
If we insert (\ref{f10}) into this equation, we get the following representation formula
\be\label{f13}
a_1(f)\left(\frac{1-p^2}{p}\right)^2\,+\,\frac{p^2}{1-p^4}
=\frac{c_0}{1-p^4}+\frac{2p-p^3}{1+p^2}c_1+\frac{p^2-p^4}{1+p^2}c_2\,=:\Phi_p(\omega).
\ee
Our aim is to prove the inequality
\be\label{f14}
|\Phi_p(\omega)|\,\leq\,\frac{1}{1-p^4},
\ee
where $\omega$ is as above. Obviously, it is sufficient to consider functions
$\omega$ holomorphic on $\overline{\ID}$. For them, we can represent the functional
$\Phi_p$ in the form
\be\label{f15}
\Phi_p(\omega)=\frac{1}{2\pi i}\int_{\partial \ID}\kappa_p(z)\omega(z)\,dz,
\ee
where
$$
\kappa_p(z)=\frac{1}{(1-p^4)(z-p)}+\frac{2p-p^3}{(1+p^2)(z-p)^2}
+\frac{p^2-p^4}{(1+p^2)(z-p)^3}.
$$
The functional $\Phi_p$ remains unchanged, if we replace in (\ref{f15}) the kernel
$\kappa_p$ by a rational function $K_p$ that has the same singular part at the point
$p$ as $\kappa_p$ and is holomorphic elsewhere in $\overline{\ID}$. Let
\beqq
K_p(z)\,&\pay=&\pay\,\frac{1}{1-p^4}\left(\frac{1}{z-p}\,+\,\frac{p}{1-pz}\right)
\,+\,\frac{2p-p^3}{1+p^2}\left(\frac{1}{(z-p)^2}\,+\,\frac{1}{(1-pz)^2}\right)\,\\
&\pay &\pay +
\frac{p^2-p^4}{1+p^2}\left(\frac{1}{(z-p)^3}\,+\,\frac{z}{(1-pz)^3}\right).
\eeqq
A lengthy but straightforward evaluation of $K_p$ on the unit circle results
in the following identity

\vspace{8pt}
$\ds e^{i\theta}K_p(e^{i\theta})(1+p^2)|1-pe^{i\theta}|^6\,$
\beqq
&\pay=&\pay (1-2p\cos\,\theta+p^2)^2\\
&\pay&\pay +(2p-p^3)(-4p+2(1+p^2)\cos\,\theta)(1-2p\cos\,\theta+p^2)\\
&\pay&\pay + (p^2-p^4)(4(\cos\,\theta)^2\,-\,(2p^3+6p)\cos\,\theta\,-2+6p^2)\\
&\pay=&\pay 4p^4(-2+p^2)(\cos\,\theta)^2\,+\,4p^3(3-p^2)\cos\,\theta\\
&\pay&\pay  + 1-8p^2+5p^4-2p^6\\
&\pay:=&\pay Q_p(\cos\,\theta),
\eeqq
where $\theta\in [0,2\pi]$. The function $Q_p(x)$ has its local maximum at the point
$$ x_M(p)\,=\,\frac{3-p^2}{2p(2-p^2)}.
$$
Since $x_M(p)>1$ for $p\in (0,1),$ we get
$$Q_p(\cos\,\theta)\geq Q_p(-1)=1-8p^2-12p^3-3p^4+4p^5+2p^6:=S(p),\quad\theta\in [0,2\pi].
$$
From $S'(p)<0$ for $p\in (0,1]$ and $S(1-\sqrt{2}/2)=0$ we conclude that
$$e^{i\theta}K_p(e^{i\theta})\geq\,0,\quad \theta\in [0,2\pi]\,\,\mbox{\rm{and}}
\,\,p\in\left(0,1-\frac{\sqrt{2}}{2}\right].
$$
Hence the desired inequality (\ref{f14}) results from the following chain of relations
\beqq
\left|\frac{1}{2\pi i}\int_{\partial \ID}K_p(z)\omega(z)   \,dz\right|
& \leq & \frac{1}{2\pi}\int_0^{2\pi}\left|e^{i\theta}K_p(e^{i\theta})\right|
\,d\theta\,\|\omega\|_{\infty}\\
&=& \frac{1}{2\pi}\int_0^{2\pi}e^{i\theta}K_p(e^{i\theta})
\,d\theta\,\|\omega\|_{\infty}\\
&\leq &\frac{1}{2\pi}\int_0^{2\pi}e^{i\theta}K_p(e^{i\theta})\,d\theta\,\\
& =& \frac{1}{2\pi i}\int_{\partial \ID}K_p(z)\,dz\\
&=& \frac{1}{1-p^4}.
\eeqq
This proves the inequality (\ref{f14}) and therefore (\ref{f12}).

For the proof that any point in the disc described by (\ref{f12}) occurs as the Laurent
coefficient $a_1(f)$ of a function $f\in Co(p)$ we may use the same functions as in the analogous
situation in the proof of Theorem \ref{th2}.

To prove the second assertion of Theorem \ref{th4} we observe that in the above chain equality is
attained everywhere if $\omega(z)\equiv 1$.
If we apply the  theory of extremum problems for linear functionals on  $H^{\infty}$
to the linear functional $\Phi_p$ (compare in particular \cite[Theorem 8.1]{D}), we see
that there is a unique extremal function $\omega_E$ such that
$$\max\{|\Phi_p(\omega)|\mid \omega \in H^{\infty}, \|\omega\|_{\infty}\leq 1\}\,
=\,\Phi_p(\omega_E).
$$
The above considerations show that in our case $\omega_E(z)\equiv 1$. This implies
that equality in (\ref{f14}) is attained if and only if  $\omega(z)\equiv e^{i\theta}$
for some $\theta \in [0,2\pi)$.
This concludes the proof of Theorem \ref{th4}.
\epf

For the remaining values of $p$ we can show that an improved version of
Theorem \ref{th3} is valid.

\bthm\label{th5}
Let $p\in(1-\frac{\sqrt{2}}{2},1)$. Then the closed convex hull of the class $Co(p)$
is a proper subset of the class of functions defined by {\rm (\ref{f3})}.
\ethm \bpf
For the proof, we use the same functions as in the proof of Theorem \ref{th3}.
For the coefficients of the Taylor expansion of $\omega_x$,
defined by (\ref{f11}), at the point $p$, we compute
$$c_0=-x,\quad c_1=-\,\frac{1-x^2}{1-p^2},\quad\mbox{\rm{and}}\quad
c_2=-\,\frac{1-x^2}{(1-p^2)^2}(p-x).
$$
If we insert these identities into (\ref{f13}), we derive the following expression for
the Laurent coefficients $a_1(f)$
$$a_1(f)\,=\,-\,\frac{p^2}{(1-p^2)^3}\left(1\,+\,\frac{1-x}{1+p^2}
\left(-1+(1+x)(2p-p^2x)\right)\right).
$$
Let $R_p(x)\,=\,-1+2p+x(2p-p^2)-p^2x^2. $
Because of
$$R_p(1)\,=\,-1+4p-2p^2\,>\,0\quad\mbox{\rm{for}}\,\,p\in\left(1-\frac{\sqrt{2}}{2},1\right),
$$
we see that there exist $x\in (0,1)$ such that for the corresponding functions
$f$ the inequalities
$$a_1(f)\,<\,-\,\frac{p^2}{(1-p^2)^3}
$$
are valid. Hence, according to Theorem C, these functions $f$ do not belong to the
closed convex hull of $Co(p)$.
\epf

\end{document}